\documentclass[conference, 10pt]{IEEEtran}

\usepackage{tikz}
\usetikzlibrary{decorations.markings}
\usetikzlibrary{arrows,calc}
\usepackage{pgfplots}
\usepackage{amscd,epsfig,bbm,graphics}

\usepackage{amsmath}
\usepackage{amsfonts}
\usepackage{amssymb}
\usepackage{titletoc}

\usepackage{stmaryrd}

\usepackage{ltugcomn}

\begin{document}

\title{A Stochastic Model for Car-Sharing Systems}

\author{\IEEEauthorblockN{Christine Fricker}
\IEEEauthorblockA{INRIA, Rocquencourt \\ France}
\and
\IEEEauthorblockN{Cedric Bourdais}
\IEEEauthorblockA{INRIA, Rocquencourt \\ Ecole Polytechnique, Palaiseau \\ France}}

\maketitle

\begin{abstract}
Vehicle-sharing systems are becoming important for urban transportation. In these systems, users arrive at a station, pick up a vehicle, use it for a while and then return it to another station of their choice. Depending on the type of system, there might be a possibility to book vehicles before picking-up and/or a parking space at the chosen arrival station. Each station has a finite capacity and cannot host more vehicles and reserved parking spaces than its capacity. We propose a stochastic model for an homogeneous car-sharing system with possibility to reserve a parking space at the arrival station when picking-up a car. We compute the performance of the system and the optimal fleet size according to a specific metric. It differs from a similar model for bike-sharing systems because of reservation that induces complexity, especially when traffic increases.
\end{abstract}

\newtheorem{Proposition}{Proposition}
\newtheorem{Theo}{Theorem}
\newtheorem{Prop}{Property}
\newtheorem{Conj}{Conjecture}
\newtheorem{Def}{Definition}
\newtheorem{Assumption}{Assumption}
\def\C{{\mathbb C}}
\def\N{{\mathbb N}}
\def\R{{\mathbb R}}
\def\Z{{\mathbb Z}}
\def\P{{\mathbb P}}
\def\Q{{\mathbb Q}}
\def\E{{\mathbb E}}
\def\mean{{\mathbb E}}
\def\P{{\mathbb P}}
\def\pr{\mathbb{P}}
\def\ind{\mathbbm{1}}
\def\cal{\mathcal}

\section{Introduction}
\label{sec:introduction}
\subsection{Context}
Over the past decade, vehicle-sharing systems have appeared as a new answer to mobility challenges, like reducing congestion, pollution or time of travel. As a consequence, the number of cities equipped with bike-sharing sytems is almost 400 (with 60,000 to 120,000 trips a day in Paris for instance) and car-sharing systems are beginning to spread as well.

Those systems can be described as follows. For a typical bike-sharing system, users simply arrive at a station, pick a bike if there is one, otherwise leave the system, use it to go to another station and leave it there, if there are spaces available. If no space is available, the users have to find a neighbouring station to return their bike.  For car-sharing systems, users have the possibility to reserve either a parking space only or both a car and a parking space in the destination station before picking-up the car. Again, they can only do this if there are cars and spaces available. Otherwise they leave the system. This basic description could be refined, especially users, instead of leaving the system unhappy, could look for a neighbouring station with available resources.  The lack of resources is one of the major issues which oblige operators to perform \textit{regulation} to maintain the reliability of the service against other transportation modes (see \cite{Bogenberger-1} and \cite{Carlier-1}).

\subsection{Related works}
Due to the success of vehicle sharing programs and the crucial problem of allocation of ressources they face, both vehicles and parking spaces, there have been many studies on redistribution in bike-sharing systems. See~\cite{Chemla} and reference therein. Much of the work concerns optimization applied to deterministic models.

Stochastic approaches have been developed recently for bike-sharing systems. For pioneer papers, see~\cite{Fayolle-7} and \cite{George-1} for a model with stations with infinite capacities, where the problem of lack of parking spaces  is avoided. In~\cite{Fricker}, the first homogeneous model addressing the problem of the availability of parking slots  aims to investigate the influence of parameters such as demand or capacity on the performance of the system, defined in terms of the proportion of \textit{problematic stations} (empty or full  stations). The paper~\cite{Fricker} deals with a set of $N$ identical stations with capacity $K$. Users arrive at rate $\lambda$ in each station, pick-up a bike if there is one, join the pool of riding users and then return their bike in a station chosen at random, after a trip duration with exponential distribution with parameter $\mu$. If the station is saturated, they reattempt with the same strategy until they succeed. The authors prove that, in this model, the  minimum stationary proportion of problematic stations as the system size gets large is $2/(K+1)$. This minimum is reached when the average number of bikes per station $s$ equals $K/2+\lambda/\mu$.

The paper~\cite{Fricker} also assesses the consequences of different incentive policies such as choosing when returning the bikes between two stations the one with  less bikes, which drastically reduces the proportion of problematic stations. In another paper~\cite{Heterogeneous}, they show how to extend the results on the basic model to a  inhomogeneous framework, considering subsets of stations with the same parameters. In those papers, the Markov process considered is the empirical distribution of the number of bikes in stations, \textit{i.e.} a vector with the proportion of stations with $k$ bikes at time $t$. The mean-field limit gives an ODE satisfied by this empirical distribution when the number of stations $N$ tends to infinity. Then,  it is proved that   the solutions of this ODE converge with time to a unique equilibrium point, and that this equilibrium point is the concentration point of the invariant measure of the empirical distribution process as the system  gets large. It gives here that, in both homogeneous and inhomogeneous cases, the  number of bikes in a station follows a truncated geometric distribution at equilibrium, when the system gets large.

In an upcoming paper~\cite{Tibi}, Tibi and Fricker prove directly the convergence of finite marginals of the invariant measure of the Markov state process, \textit{i.e.} the vector of length $N$ containing the number of bikes of each station at time $t$.  Thanks to irreducibility and reversibility, it has a product form and, using a local limit theorem, this gives moreover that the number of bikes at a fixed number of stations are independent in the limit $N \rightarrow \infty$. This property is always difficult to reach via the mean field method. An interesting result, original as far as we know,  in~\cite{Tibi} is that the Markov process including both stations and routes states has still a product-form  invariant measure, for  the model with  fixed size. It extends the result of George and Xia~\cite{George-1},  for the  scaling considered along our paper, where both numbers of stations and bikes tend to infinity, and for finite capacities.

Waserhole et al.~\cite{Jost}  propose a Markovian model for car-sharing systems  where they use a fluid limit approach to be able to maximize the average profit of the operator through revenue management techniques. 
The paper \cite{Kaspi-1} investigates the impact of reservation through different policies in vehicule-sharing systems.

\subsection{Outline of the paper}
To our knowledge, our paper presents the first stochastic analysis of a large-scale car sharing system.

We present an homogeneous model  with $M$ cars and $N$  stations with capacity $K$ where users have to reserve a parking space in a destination station chosen at random when they pick-up a car. Users arrive at rate $\lambda$ in each station. If they find a car and if there is some free space in their  destination station, they pick-up a car and make the reservation at the same time.  Otherwise, they leave the system. 

Unlike bike-sharing models, neither the empirical measure process nor the state process are reversible. Even though we know that there exists an invariant measure for a fixed $N$ for these irreducible finite-state space Markov processes, it is quite untractable. Thus the aim is to obtain large-scale asymptotics to understand the large-scale behavior of the system.  It means that $M$ and $N$ are large with the number of cars per station tending to a constant $s$. This sizing parameter is a key parameter of the system. For that, we are still able to adapt the mean-field limit approach in~\cite{Fricker} when the system gets large. The main difference is that, due to the reservation, a station is described by two components: the vehicules and the reserved parking spaces which makes the analysis of the equilibrium more tedious.  Indeed, recently,
large-system asymptotics have been successfully applied in many contexts in com-
munication systems (see 
 \cite{Simatos}, \cite{gast2010} and others). But it is, as far as we know,  the first system where the underlying process gouverning the dynamical system is not one-dimensional.

The analysis of the equilibrium point gives the following: when the system gets large, the steady-state numbers of vehicles and reserved spaces at each station have a product form distribution with a geometric and a Poisson terms on a constraint space. The two parameters involved are solutions of two fixed point equations.  This allows us to study the behavior of the system. As intuition suggests it, we prove that all stations are problematic when the average number of vehicles per station equals the capacity, which makes quite a difference with the homogeneous model for bike-sharing systems, or when there are no vehicles. We investigate the fleet sizing problem. We obtain asymptotics in the two cases of light and heavy traffic. We prove that in light traffic case, reservation has little impact on performance, unlike the heavy traffic case. The main difference with bike-sharing systems is that the best performance degrades with traffic and that the corresponding fleet size remains under the overall capacity. 

Then, a model with reservation of both resources, vehicles and parking spaces, is proposed. An approximated simple model is studied and our simulations show its relevancy. The analysis shows that the performance of such a system  behaves as in a system with  simple reservation but with more traffic. 

Section~\ref{sec:SystemModel-MeanField} describes the model, the Markovian process and the  limiting ODE. In Section~\ref{sec:SteadyStateBehaviour}, we prove the uniqueness of the equilibrium point of the ODE. In Section~\ref{sec:Performance}, we analyze the equilibrium and give performance results. Section~\ref{sec:DoubleReservation} deals with the reservation of both resources. Then Section~\ref{sec:Conclusion} concludes the paper.

\section{Stochastic Model}
\label{sec:SystemModel-MeanField}
\subsection{Main Notations}

\begin{tabular}{lp{6cm}}
$N$ & Number of stations. \\ 
$M_N$ & Total number of vehicles. \\
$s_N$ & Average number $M_N/N$ of cars per station. .\\ 
$K$ & Number of parking spaces in a station, also called capacity of the station. \\ 
$\lambda$ & Arrival rate of users at a station. \\ 
$1/\mu$ & Average trip time. \\ 
$V_i^N(t)$ & Number of vehicules in station $i$ at time $t$. \\ 
$R_i^N(t)$ & Number of parking spaces reserved in station $i$ by users travelling at time $t$. \\
$Y^N_{k,l}(t)$ & Proportion of   stations with $k$ cars and $l$ reserved parking places at time $t$. \\
$y_{k,l}(t)$ & Limit of $Y^N_{k,l}(t)$ as N tends to infinity (described by an ODE). \\
$\overline{y}$ & Equilibrium point of the corresponding ODE.\\
$\chi$ & Space of  station states.
\end{tabular}

\subsection{Model description}

The system is a set of $N$  stations with capacity $K$ with $M_N$ vehicles. As $N$ tends to infinity, $s_N = M_N/N$ tends to $s$ called the average number of vehicles per station. At a given station with capacity $K$, users are supposed to arrive with rate $\lambda$. An arriving user at station say $i$ has a destination say $j$ chosen at random.  If there is no car available in the station origin $i$ \textit{or} if there is no available parking space in the station destination $j$, the unhappy user leaves the system. Otherwise, she picks-up a car in station $i$ and simultaneously makes a reservation in station $j$. Then, the journey between station $i$ and station $j$ takes an exponentially distributed time with mean $1/\mu$. The user returns her car at station $j$ and leaves the system.

Note that unlike bike-sharing systems (\cite{Fricker}), users do not have to look for a station to return their vehicle until they find a parking space available, thanks to reservation.

In this paper, we focus on this homogeneous model. But we are able to easily extend the result to a heterogeneous model consisting in a finite number of clusters with for a station of cluster $i$ a  capacity $K_i$, an arrival rate $\lambda_i$, and probability $p_i$ of choosing a destination in cluster $i$, instead of $1/N$.

\subsection{The Empirical Measure Process}
Let us define $\chi= \{ (k,l) \in \N^2, k+l \leq K \}$ and let $Y^N_{k,l}(t)$ be the proportion of the $N$ stations with $k$ vehicles and $l$ reserved parking spaces at time $t$, that is

$$Y^N_{k,l}(t) = \frac{1}{N} \sum\limits_{i=1}^N \textbf{1}_{\{V_i^N(t)=k;R_i^N(t)=l\}}$$
where $V_i^N(t)$ is the number of vehicles and $R_i^N(t)$ the number of reserved parking spaces at station $i$ at time $t$. If $\cal{P}(\chi)$ is the set of probability measures on $\chi$, let
$$\mathcal{Y} = \{ y \in \mathcal{P}(\chi), \sum\limits_{(k,l)\in \chi} (k+l) y_{k,l} = s \}$$
and
\begin{align*}
\mathcal{Y}^N = \{ \vphantom{\sum\limits_{(k,l)\in \chi}} y \in \mathcal{P}(\chi), y_{k,l} \in   \frac{\N}{N}, 
 \sum\limits_{(k,l)\in \chi} (k+l) y_{k,l} = s_N \} \text{.}
\end{align*}
As the system is homogeneous, the process $(Y^N(t)) = (Y^N_{k,l}(t))_{k,l \in \chi}$ is a Markov process on $\mathcal{Y}^N$, described as follows.
Suppose the process $(Y^N(t))$ is at state $(y_{k,l})_{(k,l) \in \chi}$. There are two different types of transitions:
\begin{itemize}
\item \textbf{Cars picked up}. The arrival rate of users in a station in state $(k,l)$ is $\lambda y_{k,l}N$ if $k>0$. A user, who arrives at a station, makes at the same time a reservation in a station in state $(k', l')$ with probability $y_{k',l'}$ if $k'+l'<K$. Therefore, the transition rate is of $\lambda y_{k',l'}y_{k,l}N$ if $k>0, k'+l'<K$. The arrival causes $y_{k,l}$ to decrease by $1/N$ and $y_{k-1,l}$ to increase by $1/N$. And the reservation causes $y_{k',l'}$ to decrease by $1/N$ and $y_{k',l'+1}$ to increase by $1/N$.
\item \textbf{Cars returned}. When a car arrives at its reserved parking space in a station in state $(k,l)$,  $y_{k,l}$  decreases by $1/N$ and $y_{k+1,l-1}$  increases by $1/N$, as the reserved space is {\em replaced} by a car. Considering stations in state $(k,l)$, the number of reserved parking spaces is $lNy_{k,l}$. As the trips are exponentially distributed with mean $1/\mu$, this transition occurs at rate $\mu l N y_{k,l}$.
\end{itemize}

The jump matrix $Q^N$ of the process $(Y^N(t))$ is thus given by, for $y,\,y'\in \mathcal{Y}^N$,
\begin{align*}
 \begin{cases} Q^N( y , y') = \lambda y_{k',l'}y_{k,l} \;\;\;\mbox{ if }k>0,\; k'+l'<K  \\
 Q^N( y , y + \frac{1}{N}(\mathbf{e}_{k+1,l-1} - \mathbf{e}_{k,l}) )  = \mu l N y_{k,l}  \end{cases} 
\end{align*}
where $y'= y + \frac{1}{N}(\mathbf{e}_{k-1,l} - \mathbf{e}_{k,l}+\mathbf{e}_{k',l'+1} - \mathbf{e}_{k',l'})$ and $(\textbf{e}_{k,l})_{(k,l) \in \chi}$ are the vectors of the canonical basis of $\mathbb{R}^{|\chi|}$.

 It is thus irreducible and on a finite set $\mathcal{Y}^N$, which implies it admits an invariant measure denoted by $\Pi^N$. As $\Pi^N$ is analytically untractable, we use a mean-field limit when $N$ tends to $+\infty$ to give a large-scale asymptotic. 
For that, we show that $(Y^N(t))_{0 \leq t \leq T}$ tends to $y(t)_{0 \leq t \leq T}$ solution of an ODE. Then, in Section~\ref{sec:SteadyStateBehaviour}, we  prove that the ODE admits a unique equilibrium point $\overline{y}$. This equilibrium will be analyzed in Section~\ref{sec:Performance}.

\subsection{Dynamical System}
\label{sec:MeanField}
By standard arguments, for $T>0$, $(Y^N(t))_{t\in [0,T]}$ converges in distribution to $(y(t))_{t\in [0,T]}$ unique solution with $y(0)$ fixed of the following ODE
\begin{multline}\label{age}
\dot{y}(t) = \sum\limits_{(k,l) \in \chi} y_{k,l}(t) \left(\sum\limits_{(k',l') \in \chi} y_{k,l}(t) \mu l (\textbf{e}_{k+1,l-1} - \textbf{e}_{k,l}) \right.
		\\  \hspace{-5mm} \left. +\lambda y_{k',l'} \mathbf{1}_{k>0;k'+l'<K} (\mathbf{e}_{k-1,l}-\mathbf{e}_{k,l}+\mathbf{e}_{k',l'+1}-\mathbf{e}_{k',l'}) \vphantom{\sum\limits_{(k,l) \in \chi}} \right).
\end{multline}

The first term corresponds to the rate at which users return cars at the reserved parking space, and the second term corresponds  to the rate of simultaneous arrival and reservation. Let us introduce some additional notations. Let
\begin{align}\label{defs}
y_S(t) = \sum\limits_{k+l=K} y_{k,l}(t), \; y_{k,.}(t) = \sum\limits_{l=0}^{K-k} y_{k,l}(t)
\end{align}
 be respectively the limiting proportion of saturated stations at time $t$ and
the limiting proportion of stations with $k$ vehicles at time $t$ ($0\leq k\leq K$).
Then, by splitting the second term of its right-hand side, the ODE~\eqref{age} rewrites
\begin{multline*}
 \dot{y}(t) =  \sum\limits_{(k,l) \in \chi} y_{k,l}(t)\left(  \lambda (1-y_S)(\mathbf{e}_{k-1,l}-\mathbf{e}_{k,l}) \mathbf{1}_{k>0} \right.
\\\left. + \lambda (1-y_{0,.}) (\mathbf{e}_{k,l+1}-\mathbf{e}_{k,l}) \mathbf{1}_{k+l<K} +  \mu l (\textbf{e}_{k+1,l-1} - \textbf{e}_{k,l}) \right).
\end{multline*}
of the form
\begin{align}\label{ODEModel1}
\dot{y}(t) = y(t)L_{y(t)}
\end{align}
 where $ y(t)L_{y(t)}$ is the product of the vector $y(t)$ by the jump matrix $L_{y(t)}$ defined in the following section. This means that, when $N$ tends to infinity, the empirical distribution $y(t)$ of the stations evolves in time as the distribution of some non-homogeneous Markov process on $\chi$, whose jumps are given by $L_{y(t)}$, updated by the current distribution $y(t)$.

\section{Steady-State Behavior}
\label{sec:SteadyStateBehaviour}
\subsection{Probabilistic Interpretation of the ODE}
\label{sec:InterpretationODE}

\begin{figure}[h]
\centering
\begin{tikzpicture}[scale=0.8]
\draw[->] (1.8,0) -- (3.2,0);
\draw[->] (4.8,0) -- (6.2,0);
\draw[->] (9.8,0) -- (11.2,0);

\draw[-] (4.5,-1.5) -- (4.5,1.5);
\draw[-] (3.5,-1.5) -- (3.5,1.5);
\draw[-] (3.5,1.5) -- (4.5, 1.5);
\draw[-] (3.5,1) -- (4.5, 1);
\draw[-] (3.5,.5) -- (4.5,.5);
\draw[-] (3.5,0) -- (4.5,0);

\draw[-] (6.5,.5) -- (9.5,.5);
\draw[-] (6.5,-.5)-- (9.5,-.5);
\draw[-] (9.5,.5) -- (9.5, -.5);
\draw[-] (9,.5) -- (9, -.5);
\draw[-] (8.5,.5) -- (8.5, -.5);
\draw[-] (8,.5) -- (8, -.5);

\node (arrival) at (2.5,0.5) {$\lambda(1-y_{0,.})$};
\draw (arrival);
\node (middle) at (4.8,1.5) {$\mu$};
\draw (middle);
\node (out) at (10.3,0.8) {$\lambda (1-y_S)$};
\draw (out);

\node (Reservations) at (4,2) {R: $M/M/\infty$};
\draw (Reservations);

\node (Vehicules) at (8,1) {V: $M/M/1$};
\draw (Vehicules);

\end{tikzpicture}
\caption{A typical station as a tandem of two queues with overall capacity $K$.}
\label{fig:MMqueuesModel1}
\end{figure}
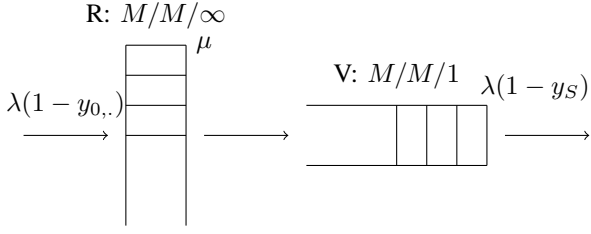

The jumps of the Markov process with generator $L_y$ on $\chi= \{ (k,l) \in \mathbb{N}^2, k+l \leq K \}$ defined in Section~\ref{sec:SystemModel-MeanField} can be seen as those of the number of customers in two coupled queues (see Figure~\ref{fig:MMqueuesModel1}): a typical station can be described as a tandem of two queues: a $M/M/\infty$ queue for reservations, with arrival rate $\lambda(1-y_{0,.})$ and service rate $\mu$, and a $M/M/1$ queue for vehicles, where customers come from the former queue, with service rate $\lambda (1-y_S)$. Moreover, this system is a loss system: the total number of customers is less than or equal to $K$. Then the problem is now to find the stationary measure of this system of two queues in tandem.

Let $\rho_R$ be the arrival-to-service rate ratio of the $M/M/\infty$ queue, and $\rho_V$ of the $M/M/1$ queue, defined as
\begin{eqnarray}
\rho_R(y) &=& \frac{\lambda}{\mu}(1-y_{0,.}), \label{eq:RhoRy}
\\  \rho_V(y) &=& \frac{1-y_{0,.}}{1-y_S}. \label{eq:RhoVy}
\end{eqnarray}

In this case, it is well-known (see~\cite{Robert} for example) that the invariant probability measure  $\pi(y)$ associated to $L_y$ has a product form, given for $(k,l) \in \chi$ by
\begin{multline}
\pi_{k,l}(\rho_V(y), \rho_R(y)) = \frac{1}{Z(\rho_V(y), \rho_R(y))}\frac{\rho_R(y)^l}{l!}\rho_V(y)^k  \label{MeasureModel1}
\end{multline}
where $Z(\rho_V(y), \rho_R(y))$ is such that 
$$\sum\limits_{(k,l) \in \chi} \pi_{k,l}(\rho_V(y), \rho_R(y)) = 1.$$

With a slight abuse of notations, we will use both $\pi(y)$ and $\pi(\rho_V(y), \rho_R(y))$ in the following.

Because $\pi(y) \in \mathcal{Y}$, we also have
\begin{equation}
\label{eq:sMoyenneY}
s = \sum\limits_{(k,l) \in \chi} (k+l)\pi_{k,l}(y)  \text{.}
\end{equation}
which is equivalent to
$s =  \mathbb{E}(V+R)$
where $(R,V)$ is a random variable with distribution $\pi(y)$.

To investigate the steady-state behaviour of the model, we study the equilibrium points $\overline{y}$ of the ODE~\eqref{ODEModel1}. Such equilibrium points $\overline{y}$ are the $y$ satisfying
\begin{equation}
\label{eq:fixedPoint}
\pi(\rho_V(y), \rho_R(y)) = y
\end{equation}
and equations~\eqref{eq:RhoRy},~\eqref{eq:RhoVy} and~\eqref{eq:sMoyenneY}.
Therefore, it is also equivalent to finding couples $(\rho_V, \rho_R)$ such that $\bar{y}=\pi(\rho_V, \rho_R)$ satisfying
\begin{eqnarray}
\rho_R &=& \frac{\lambda}{\mu}(1-\pi_{0,.}(\rho_V, \rho_R)) \label{eq:RhoR}
\\  \rho_V &=& \frac{1-\pi_{0,.}(\rho_V, \rho_R)}{1-\pi_S(\rho_V, \rho_R)} \text{.} \label{eq:RhoV}
\\ 	s 		&=& \mathbb{E}(R+V) \label{eq:sMoyenne}
\end{eqnarray}
where $(R,V)$ is a random variable with distribution $\pi(\rho_V, \rho_R)$.

In conclusion, the question of finding a measure on $\chi$ is reduced to find a couple of real numbers $(\rho_V, \rho_R)$. The two following sections prove the existence and uniqueness of the equilibrium point.


\subsection{Existence of an Equilibrium Point}
\label{sec:ExistenceFixedPoint}

Due to the  Section~\ref{sec:InterpretationODE}, we know that finding an equilibrium point $\overline{y}$ of the ODE~\eqref{ODEModel1} is equivalent to finding  $\overline{y}$ solution of equation~\eqref{eq:fixedPoint}.
Since $s < K$,
for each $y \in \mathcal{Y}, \ y_S <  1$.
Therefore, by definition (see equations~\eqref{eq:RhoRy}, ~\eqref{eq:RhoVy} and~\eqref{MeasureModel1}), $y \mapsto (\rho_V(y),\rho_R(y))$ is  continuous from the compact and convex set $\mathcal{Y}$ to itself.
 Thus, the existence of $\bar{y}$ is given by Brouwer's theorem.


\subsection{Uniqueness of the Equilibrium Point}
\label{sec:UniquenessFixedPoint}

Compared to the bike-sharing model in~\cite{Fricker}, reservation makes calculations much more tedious. But we still prove the uniqueness of the equilibrium point. Let us  present the main steps of the proof.

 First, we prove that Equation~\eqref{eq:RhoV} is  true for  any $(\rho_V,\rho_R)$. Second, we consider equation~\eqref{eq:RhoR} as an implicit function between $\rho_R$ and $\rho_V$ which gives $\rho_V$ as a strictly increasing function of $\rho_R$. The third step consists in studying the monotony of the right-hand side of equation~\eqref{eq:sMoyenne} as a function of $\rho_V$ and $\rho_R$ (as for bike-sharing systems, it increases), and then, in using the diffeomorphism between $\rho_R$ and $\rho_V$ to conclude on the uniqueness of the solution of~\eqref{eq:RhoR} 
and \eqref{eq:sMoyenne}.

From now, we assume~\eqref{eq:fixedPoint} to be satisfied since we know there exists solutions. We need to find the solutions $(\rho_V, \rho_R)$ of equations~\eqref{eq:RhoR},~\eqref{eq:RhoV} and~\eqref{eq:sMoyenne}.  Let us first notice that necessarily $(\rho_V, \rho_R) \in \Gamma = [0, \infty \mathclose{[} \times [0, \lambda/\mu \mathclose{[}$ according to equations~\eqref{eq:RhoR} and~\eqref{eq:RhoV}.

The partition function $Z(\rho_V, \rho_R)$ of the invariant probability measure  plays an important part in the following, so we   mention some of its properties.

\begin{Prop}
\label{PropZ}
(i) For any $(\rho_V, \rho_R) \in \Gamma=[0, \infty \mathclose{[} \times [0, \lambda/\mu \mathclose{[}$,
\begin{align*}
Z(\rho_V, \rho_R) &= \sum\limits_{(k,l) \in \chi} \rho_V^k \frac{\rho_R^l}{l!}  \text{,}\\
 \frac{\partial Z}{\partial \rho_R} &= (1-\pi_S)Z  \text{,}\\
Z \left(\frac{1}{\rho_V}, \rho_R \right) \rho_V^K &= Z(\rho_V,\rho_R\rho_V)  \text{ if } \rho_V>0.
\end{align*}
(ii) For any $(\rho_V, \rho_R)$ solution of equation~\eqref{eq:RhoR}, 
$$\left(\frac{\lambda}{\mu}-\rho_R \right)Z(\rho_R, \rho_V) = \frac{\lambda}{\mu}\sum\limits_{l=0}^K \frac{\rho_R^l}{l!}  \text{.}$$
\end{Prop}
\begin{IEEEproof}
Simple algebra gives the result. For the third equation of $(i)$, using the change of indexes $k'=K-k-l$,
\begin{align*}
Z \left(\frac{1}{\rho_V}, \rho_R \right) \rho_V^K &= \sum\limits_{(k,l) \in \chi} \rho_V^{K-k} \frac{\rho_R^{l}}{l!} = \sum_{l=0}^K \frac{\rho_R^{l}}{l!}\sum_{k=0}^{K-l}\rho_V^{K-k}  \\
&= \sum_{l=0}^K \frac{\rho_R^{l}}{l!}\sum_{k'=0}^{K-l}\rho_V^{k'+l} =\sum\limits_{(k',l) \in \chi}\frac{\rho_R^{l}}{l!}\rho_V^{k'+l} \\			
	&= Z(\rho_V, \rho_R\rho_V).
\end{align*}
\end{IEEEproof}

\begin{Prop}
\label{Prop:RhoV1}
Equation~\eqref{eq:RhoV} is true for all $(\rho_V, \rho_R) \in \Gamma=[0, \infty \mathclose{[} \times [0, \lambda/\mu \mathclose{[}$.
\end{Prop}
\begin{IEEEproof}
Using definitions given by~\eqref{defs},
\begin{align*}
\rho_V (1-\pi_S(\rho_V, \rho_R)) 	&= \frac{1}{Z(\rho_V, \rho_R)}\rho_V\sum\limits_{k+l < K} \frac{\rho_R^l}{l!}\rho_V^k \\
&= \frac{1}{Z(\rho_V, \rho_R)}\sum\limits_{k>0, k+l \leq K} \frac{\rho_R^l}{l!}\rho_V^k \\
							&= 1-\pi_{0,.}(\rho_V, \rho_R).
\end{align*}\end{IEEEproof}

Property~\ref{Prop:RhoV1} ends the first step of the proof: It allows us to focus on equations~\eqref{eq:RhoR} and~\eqref{eq:sMoyenne} only.  It means that the set of equations~\eqref{eq:RhoR}, ~\eqref{eq:RhoV} and ~\eqref{eq:sMoyenne} is equivalent to the set of equations ~\eqref{eq:RhoR} and ~\eqref{eq:sMoyenne} .

\begin{Theo}
\label{Theorem:diffeomorphism}
There exists a strictly increasing diffeomorphism $\phi: \mathopen{[}0, \lambda/\mu \mathclose{[} \rightarrow [0, \infty \mathclose{[}$, and $\psi=\phi^{-1}$, such that  $(\rho_V, \rho_R)$ is a solution of equation~\eqref{eq:RhoR} if and only if $\rho_V = \phi(\rho_R)$.
\end{Theo}
\begin{IEEEproof}
It uses the {\em implicit function theorem}. Due to Property~\ref{PropZ} $(ii)$, $(\rho_V, \rho_R) \in \Gamma$  is solution of equation~\eqref{eq:RhoR} if and only if $f(\rho_V, \rho_R) = 0$ where $f$ is the  $\mathcal{C}^{\infty}$ function defined by
\begin{align}\label{f}
f(\rho_V, \rho_R) = \left(\frac{\lambda}{\mu}-\rho_R \right)Z(\rho_R, \rho_V) - \frac{\lambda}{\mu}\sum\limits_{l=0}^K \frac{\rho_R^l}{l!}.
\end{align}

And we already know by Section~\ref{sec:ExistenceFixedPoint} that there exists a solution of equation~\eqref{eq:RhoR}. Let $(\rho_V, \rho_R)$ be such a solution. To prove the existence of a strictly increasing diffeomorphism, we have to show that
$$\frac{\partial f}{\partial \rho_R}(\rho_V, \rho_R) \frac{\partial f}{\partial \rho_V}(\rho_V, \rho_R) < 0 \text{.}$$

First, using Property~\ref{PropZ} $(i)$,
\begin{multline}\label{df/dr_R}
\frac{\partial f}{\partial \rho_R}(\rho_V, \rho_R) = -Z(\rho_V, \rho_R) 
\\+ \left(\frac{\lambda}{\mu}-\rho_R \right) \sum\limits_{k+l \leq K-1} \rho_V^k \frac{\rho_R^l}{l!} -\frac{\lambda}{\mu} \sum\limits_{l=0}^{K-1}\frac{\rho_R^l}{l!}
\end{multline}
thus, substracting equation~\eqref{df/dr_R} to equation~\eqref{f},
\begin{multline}\label{diff}
f(\rho_V, \rho_R)-\frac{\partial f}{\partial \rho_R}(\rho_V, \rho_R) = Z(\rho_V, \rho_R) 
\\+ \left(\frac{\lambda}{\mu}-\rho_R \right) \sum\limits_{k=0}^K \frac{\rho_V^k \rho_R^{K-k}}{(K-k)!} - \frac{\lambda}{\mu} \frac{\rho_R^K}{K!}.
\end{multline}

From equation~\eqref{f}, while $f(\rho_V, \rho_R)=0$, it holds that
\begin{align*}
\frac{\lambda}{\mu}-\rho_R=\frac{\lambda}{\mu Z(\rho_R, \rho_V)}\sum\limits_{l=0}^K \frac{\rho_R^l}{l!}.
\end{align*}
Thus, the second term of the right-hand side of equation~\eqref{df/dr_R} can be rewritten
\begin{align}\label{Majoration}
& \left(\frac{\lambda}{\mu}-\rho_R \right) \sum\limits_{k=0}^K \frac{\rho_V^K \rho_R^{K-k}}{(K-k)!}  \\
	& = \frac{\lambda}{\mu} \frac{1}{Z(\rho_V, \rho_R)} \sum\limits_{l=0}^K \frac{\rho_R^l}{l!} \sum\limits_{k=0}^K \frac{\rho_V^k \rho_R^{K-k}}{(K-k)!}  \nonumber \\
	& = \frac{\lambda}{\mu}  \frac{\rho_R^K}{Z(\rho_V, \rho_R)}  \sum\limits_{  i,j=0}^K \frac{\rho_R^{i-j}\rho_V^j}{i!(K-j)!}\nonumber \\ 
	& = \frac{\lambda}{\mu}  \frac{\rho_R^K}{Z(\rho_V, \rho_R)} \left( \sum\limits_{0 \leq j \leq i \leq K} \frac{\rho_R^{i-j}\rho_V^j}{i!(K-j)!} + \sum\limits_{j>i} \frac{\rho_R^{i-j}\rho_V^j}{i!(K-j)!} \right) \nonumber
\end{align}

For the first term in the sum of the right-hand side of equation~\eqref{Majoration}, 
\begin{align*}
\hspace{-5mm} \sum\limits_{0 \leq j \leq i \leq K} \frac{\rho_R^{i-j}\rho_V^j}{i!(K-j)!} = \frac{1}{K!} \sum\limits_{(k,l) \in \chi} \frac{K!\rho_R^{l}\rho_V^k}{(k+l)!(K-k)!}  
  \geq \frac{Z(\rho_V, \rho_R)}{K!}
\end{align*}
using that, for any $(k,l) \in \chi$, 
\begin{align*}
\frac{K!}{(k+l)!(K-k)!} = \frac{1}{l!}\prod_{i=0}^{k-1}\frac{K-i}{k+l-i} \geq \frac{1}{l!}.
\end{align*}
 Moreover, the second term in the sum  on the right-hand side of~\eqref{Majoration} is positive. Thus equation~\eqref{Majoration} yields that
\begin{equation}
\label{Majoration2}
\left(\frac{\lambda}{\mu}-\rho_R \right) \sum\limits_{k=0}^K \frac{\rho_V^K \rho_R^{K-k}}{(K-k)!} \geq  \frac{\lambda}{\mu} \frac{\rho_R^K}{K!} \text{.}
\end{equation}

Additionally, using that $Z>0$, we can deduce from equation~\eqref{diff} that $f(\rho_V, \rho_R)- \frac{\partial f}{\partial \rho_R}(\rho_V, \rho_R) > 0$. Therefore, while $f(\rho_V, \rho_R)=0$,
$$\frac{\partial f}{\partial \rho_R} (\rho_V, \rho_R)<0  \text{.}$$

It is easy to see that $\frac{\partial f}{\partial \rho_V}(\rho_V,\rho_R)$ is positive: Since $\rho_R < \lambda/\mu$ and $Z$ strictly increasing in $\rho_V$,
$$\frac{\partial f}{\partial \rho_V}(\rho_V, \rho_R) = \left(\frac{\lambda}{\mu}-\rho_R \right) \frac{\partial Z}{\partial \rho_V}(\rho_V, \rho_R)  >0  \text{.}$$

One can easily check that if $\rho_V = 0$ then $\rho_R=0$. Similarly, when $\rho_V$ tends to infinity, $\rho_R$  tends to $\lambda/\mu$ to keep $f(\rho_V, \rho_R) = 0$. It ends the proof.
\end{IEEEproof}

Theorem~\ref{Theorem:diffeomorphism}  means that the solutions of equation~\eqref{eq:RhoR} can be expressed with only one parameter. It will also be useful for a further study of equilibrium as it allows to differentiate expressions. This concludes the second step of the proof.









Once we have used equations~\eqref{eq:RhoR} and~\eqref{eq:RhoV}, let us focus now on  equation~\eqref{eq:sMoyenne}. We prove that $\E(R+V)$ as a function of $\rho_V$ and $\rho_R$ is strictly increasing in each parameter, and then using Theorem~\ref{Theorem:diffeomorphism}, we conclude to the uniqueness of the equilibrium point.

\begin{Prop}
\label{sCroissante}
The average number $\tilde{s}(\rho_V ,\rho_R)= \mathbb{E}_{}(R+V)$ of vehicles and reserved places per station,  
 where $(R,V)$ is a random variable with distribution $\pi(\rho_V ,\rho_R)$,
is a strictly increasing function of both $\rho_V$ and $\rho_R$. 
\end{Prop}
\begin{IEEEproof}
See Appendix~\ref{sec:Appendices}
\end{IEEEproof}

The following result concludes the  proof of the uniqueness of the equilibrium point  for the ODE.

\begin{Theo}
\label{Theorem:sCroissante}
For any $s>0$, there exists  a unique  $(\rho_V,  \rho_R)$ solution of both equations~\eqref{eq:RhoR} and ~\eqref{eq:sMoyenne}.

\end{Theo}
\begin{IEEEproof}
We have to prove that the function from  $[0, +\infty \mathclose{[}$ to $[0,K \mathclose{[}$ which maps
$\rho_V$ to $s_K(\rho_V, \psi(\rho_V))$ is a strictly increasing diffeomorphism.
It is true, using both Theorem~\ref{Theorem:diffeomorphism} and Property~\ref{sCroissante}, and that
$$s_K' =  \frac{\partial s_K}{\partial \rho_V} + \frac{\partial s_K}{\partial \rho_R} \psi'. $$
It ends the proof.
\end{IEEEproof}

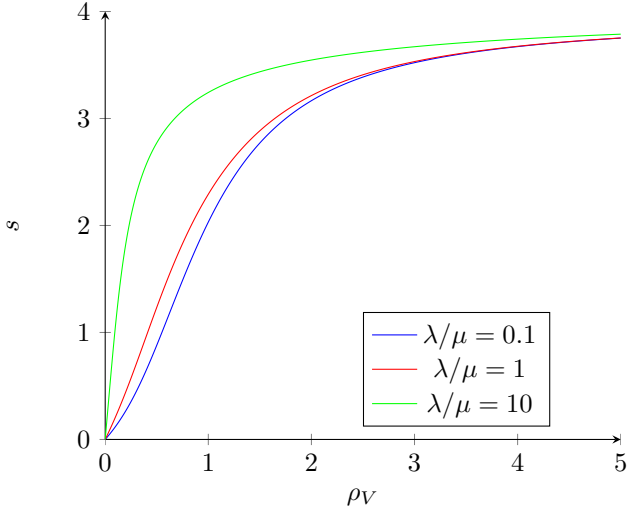
\begin{figure}
\centering

\begin{tikzpicture}[scale=1]
\begin{axis}[
title={},
area legend,
legend style={at={(0.5,0.03)},
anchor=south west},
axis x line=bottom,
axis y line=left,
xlabel={$\rho_V$},
ylabel={$s$},
xmin=0,
ymin=0,
xmax=5,
ymax=4
]

\legend{$\lambda/\mu=0.1$, $\lambda/\mu=1$, $\lambda/\mu=10$}
\addplot[blue, line legend] table {./s_fonction_rhoV_K=4_a=0.1.txt};
\addplot[red, line legend] table {./s_fonction_rhoV_K=4_a=1.txt};
\addplot[green, line legend] table {./s_fonction_rhoV_K=4_a=10.txt};
\end{axis}
\end{tikzpicture}

\caption{Diffeomorphism between $\rho_V$ and the average number of vehicles $s$ with $K=4$. $s$ is strictly increasing with traffic even though the stations' capacity is finite.}
\label{fig:sMoyenne}
\end{figure}

Figure~\ref{fig:sMoyenne} represents the diffeomorphism between $s$ and $\rho_V$ for different values of traffic $\lambda/\mu$. As we could have expected, it is globally increasing with traffic.

The proof of the uniqueness of the equilibrium point of the ODE is much more tedious than in the model for  bike-sharing systems in~\cite{Fricker} since we have two parameters related by just an implicit function.
\\\indent
This equilibrium point  also allows us to study the performance of the system in Section~\ref{sec:Performance}. We can parametrize all quantities with one parameter only and also differentiate them.


\subsection{Convergence to the Equilibrium Point}
\label{sec:ConvergenceODE}

Both the convergence of the dynamical system to its unique equilibrium point and the convergence to the steady-state empirical measure process  to the same point are out of the scope of the paper. These questions are still open.
In the bike-sharing model of \cite{Fricker}, thanks to reversibility, a Lyapunov function is found. Here, the stochastic model for car-sharing systems is much more complex: The state process of a tandem of queues is not reversible.
Recall also that Tibi \cite{Tibi} proves the convergence of  the invariant measure of the state process by its  explicit product-form. It uses a local limit theorem. For the car-sharing  model, no closed-form expression of the invariant measure can be expected.

We thus study $\overline{y}$ as the heuristic limit of both the dynamical system and the stationary empirical measure $Y^N(\infty)$. It gives the limiting stationary behavior of our model.

\section{System Performance}
\label{sec:Performance}
To manage the system, it  would be interesting to know how to choose the key parameters $K$ and $s$ to satisfy most users of the system, or to maximize profit. In this section, we discuss the choice of a metric to assess the homogeneous system performance and then give results or conjectures about optimal behavior deriving from the main  metrics. 

\subsection{Choice of the metric}
By Section~\ref{sec:SteadyStateBehaviour}, admitting the convergence as  the system size gets large of the invariant measures,  the limiting stationary distribution of the joint numbers of vehicles and reserved places in a station is $\bar{y}$. We mainly focus in our paper on a quality of service indicator, called the \textit{limiting stationary proportion of problematic stations}.

\begin{Def}[Problematic Stations]
\label{def:ProblematicStations}
Let $\overline{y}$ be the unique equilibrium point of ODE~\eqref{ODEModel1}. The stations with either no car or no parking space available are called \textit{problematic}. The limiting stationary proportion $P_b$ of problematic stations  is given by
$$P_b = \mathbb{P}(V=0 \text{ or } V+R=K)= \overline{y}_{0,.} + \overline{y}_S - \overline{y}_{0,K}$$ 
where $(V,R)$ is a random variable with distribution $\bar{y}$.
For simplification purposes, we also define $P_b^{+}$ by
$$P_b^{+}=\overline{y}_{0,.} + \overline{y}_S.$$
\end{Def}

 Note that $ P_b^+\geq  P_b $. We can also expect that $P_b^+$ is a good approximation of $  P_b$. This metric gives a convenient way to compare the performance with the bike-sharing system performance described in~\cite{Fricker}. It is also close to the proportion of unsatisfied users. Indeed, a user is satisfied if she can pick up a car at her chosen source station and can return it at her chosen destination station. In a homogeneous system, and if stations are chosen uniformly at random, this occurs with probability $(1-\overline{y}_{0,.})(1-\overline{y}_{S})$. Hence, the probability of unsatisfaction is $U=\overline{y}_{0,.} + \overline{y}_{S} - \overline{y}_{0,.}\overline{y}_{S}$. Generally, the third term is small compared to the others. In our case, we can note that $P_b \geq U$ because of  equation~\eqref{Majoration2}, key argument in the proof of Theorem~\ref{Theorem:diffeomorphism} which rewrites, using equation~\eqref{Majoration},
$$\frac{1}{Z(\rho_V,\rho_R)} \sum\limits_{l=0}^K \frac{\rho_R^l}{l!} \sum\limits_{k=0}^K \frac{\rho_V^k \rho_R^{K-k}}{(K-k)!} \geq \frac{\rho_R^K}{K!}$$
which gives, by dividing by $Z(\rho_V,\rho_R)$,
$\overline{y}_{0,.}\overline{y}_{S} \geq \overline{y}_{0,K} \text{.}$

The following equation, straightforward consequence of equations~\eqref{eq:RhoR} and~\eqref{eq:RhoV}, will be useful.

\begin{equation}
\label{eq:Pbexpression}
P_b^{+}=2 - \frac{\mu \rho_R}{\lambda} - \frac{\mu \rho_R}{\lambda \rho_V} \text{.}
\end{equation}

Therefore, the \textit{limiting proportion of problematic stations} $P_b$ is close (and an upper bound) to the \textit{limiting proportion of unsatisfied users} who can not find a car or return their car in the desired stations, which is a frequently used performance metric, and $P_b^{+}$ is a useful upper bound for both of them. It is summarized in the following proposition.
\begin{Proposition}[Metric comparison]
\begin{align*}
P^+_b\geq P_b \geq U.
\end{align*}
\end{Proposition}

One might also be interested in the number of trips achieved per unit of time. Note that  the operator revenue is more or less proportional to this quantity.

\begin{Def}[Number of Successful Trips]
Let $\overline{y}$ be the unique equilibrium point of ODE~\eqref{ODEModel1}. The limiting stationary number of successful trips per unit of time is denoted by $T$ and given by
$$T = \lambda (1-U)$$
and moreover, $T  \geq \lambda (1-P_b)$ and approximated by $\lambda (1-P_b)$.
\end{Def}



Note that in this model, the average sojourn time does not depend on the state of the stations since the parking space is reserved at the destination. Thus, we will not consider any metric related to this quantity unlike what is done for bike-sharing systems in~\cite{Fricker}.

In the following, a value of parameter $X$ referred as optimal must be understood in the sense of minimizing the proportion of problematic stations for a given traffic. Such a value is denoted by $X\text{*}$.


\subsection{Optimal Fleet Size}
We study in this section the influence of parameters on performance, especially the fleet size given by parameter $s$.

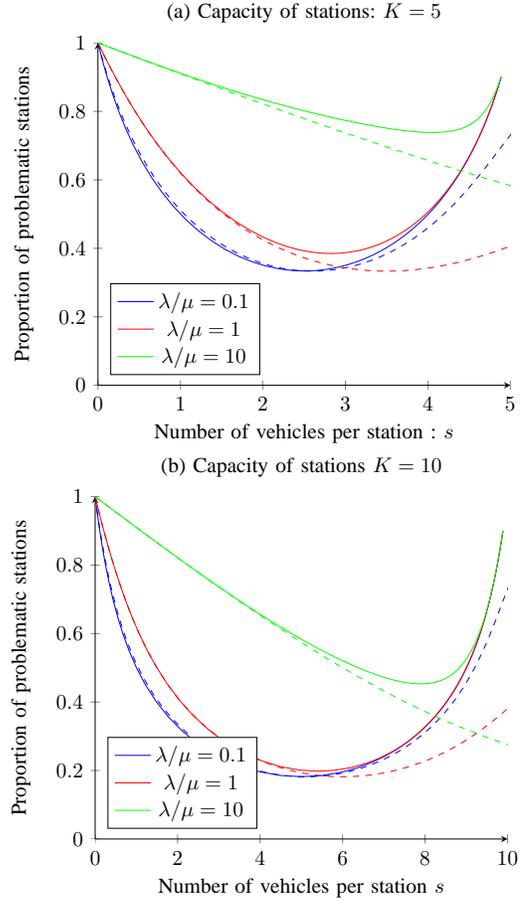
\begin{figure}[h]
\centering
\begin{tikzpicture}[scale=0.8]
\begin{axis}[
title={(a) Capacity of stations: $K=5$},
area legend,
legend style={at={(0.03,0.03)},
anchor=south west},
axis x line=bottom,
axis y line=left,
xlabel={Number of vehicles per station : $s$},
ylabel={Proportion of problematic stations},
xmin=0,
ymin=0,
xmax=5,
ymax=1
]

\legend{$\lambda/\mu=0.1$, $\lambda/\mu=1$, $\lambda/\mu=10$}
\addplot[blue, line legend] table {./Rapport_Curves_Model1_a=0.1_K=5.txt};
\addplot[red, line legend] table {./Rapport_Curves_Model1_a=1_K=5.txt};
\addplot[green, line legend] table {./Rapport_Curves_Model1_a=10_K=5.txt};

\addplot[blue, dashed, line legend] table {./Velib_K=5_a=0.1_Pbatic_stations.txt};
\addplot[red, dashed, line legend] table {./Velib_K=5_a=1_Pbatic_stations.txt};
\addplot[green, dashed, line legend] table {./Velib_K=5_a=10_Pbatic_stations.txt};
\end{axis}
\end{tikzpicture}
\hskip 10pt
\begin{tikzpicture}[scale=0.8]
\begin{axis}[
title={(b) Capacity of stations $K=10$},
area legend,
legend style={at={(0.03,0.03)},
anchor=south west},
axis x line=bottom,
axis y line=left,
xlabel={Number of vehicles per station $s$},
ylabel={Proportion of problematic stations},
xmin=0,
ymin=0,
xmax=10,
ymax=1
]

\legend{$\lambda/\mu=0.1$, $\lambda/\mu=1$, $\lambda/\mu=10$}
\addplot[blue, line legend] table {./Rapport_Curves_Model1_a=0.1_K=10.txt};
\addplot[red, line legend] table {./Rapport_Curves_Model1_a=1_K=10.txt};
\addplot[green, line legend] table {./Rapport_Curves_Model1_a=10_K=10.txt};

\addplot[blue, dashed, line legend] table {./Velib_K=10_a=0.1_Pbatic_stations.txt};
\addplot[red, dashed, line legend] table {./Velib_K=10_a=1_Pbatic_stations.txt};
\addplot[green, dashed, line legend] table {./Velib_K=10_a=10_Pbatic_stations.txt};
\end{axis}

\end{tikzpicture}
\caption{Parametric curve $\rho_V \mapsto (s, P_b)$ for $K=5$ \textbf{(a)}  and $K=10$ \textbf{(b)}. \textit{Solid lines} represent the model with reservation (car-sharing) while \textit{dashed lines} represent the model without reservation (bike-sharing).
}
\label{ParametricPb}
\end{figure}

Due to the results of the previous section, mainly Theorem~\ref{Theorem:diffeomorphism} giving the existence of $\psi$, $\bar{y}$ function of $(\rho_V,\rho_R)$, and other derived functions as performance metrics, can be   considered  as a functions of $\rho_V$ only, because $\bar{y}$ is solution of equation \eqref{eq:RhoR}. With abuse of notation, these variables will be denoted by $\bar{y}(\rho_V)$, $P_b(\rho_V)$ and $s(\rho_V)$.
Thus the proportion of problematic stations, as a function of the fleet size, is given by the parametric curve
$$\rho_V \mapsto (s(\rho_V),P_b(\rho_V)).$$
This curve will play an important role in the following, giving the behavior of the system according to the fleet size parameter. It is numerically plotted in Figure~\ref{ParametricPb}. We emphasized that this curve is not obtained by simulations. The performance metric chosen here is compared to the same for bike sharing systems (see \cite{Fricker} for details).

Further results will rely on an interesting property of symmetry  as follows. 

\begin{Prop}[Symmetry]
\label{Prop:symmetry}
For all $\rho_V > 0$,
\begin{align}\label{symP_b}
\psi(1/\rho_V) &= \frac{1}{\rho_V} \psi(\rho_V),\nonumber\\
\overline{y}_S (1/\rho_V) &= \overline{y}_{0,.} (\rho_V),\nonumber\\
P_b(1/\rho_V) &= P_b(\rho_V)\text{ and }  P_b^{+}(1/\rho_V) = P_b^{+}(\rho_V).
\end{align}
\end{Prop}
\begin{IEEEproof}
Let $\rho_V > 0$ and $\rho_R = \psi(\rho_V)$. We also define $\tilde{\rho}_R = \psi(1/\rho_V)$. Then,
\begin{align*}
\overline{y}_S(1/\rho_V) &= \frac{1}{Z(1/\rho_V, \tilde{\rho}_R)}\frac{1}{\rho_V^K} \sum\limits_{k=0}^K \frac{(\rho_V\tilde{\rho}_R)^{K-k}}{(K-k)!} \text{.}\nonumber
\end{align*}
But, with Property~\ref{PropZ}, $Z(1/\rho_V, \tilde{\rho}_R)\rho_V^K = Z(\rho_V, \rho_V \tilde{\rho}_R)$. Hence,
$$\overline{y}_S(1/\rho_V) = \frac{1}{ Z(\rho_V, \rho_V \tilde{\rho}_R) }  \sum\limits_{l=0}^K \frac{(\rho_V\tilde{\rho}_R)^{l}}{l!} \text{.}$$

Using Theorem~\ref{Theorem:diffeomorphism}, since $(\rho_V, \rho_V \tilde{\rho}_R)$ is solution of Equation~\eqref{eq:RhoR}, it holds that $\rho_V \tilde{\rho}_R = \rho_R$, which is equivalent to $\psi(1/\rho_V) = 1/\rho_V \psi(\rho_V)$. Thus $\overline{y}_S (1/\rho_V) = \overline{y}_{0,.} (\rho_V)$.

Then, with the results above,
\begin{align*}
P_b\left(\frac{1}{\rho_V}\right) &= \overline{y}_S\left(\frac{1}{\rho_V}\right) + \overline{y}_{0,.}\left(\frac{1}{\rho_V}\right)-\frac{\rho_R(\rho_V)^K}{\rho_V^K Z(\frac{1}{\rho_V}, \rho_R(\frac{1}{\rho_V}))K!}\\
 &= \overline{y}_{0,.}(\rho_V) + \overline{y}_S(\rho_V) -\frac{\rho_R(\rho_V)^K}{Z(\rho_V, \rho_R)K!}.
\end{align*}
Thus, $P_b(1/\rho_V) = P_b(\rho_V) \text{.}$
The same is true for $P_b^{+}$.
\end{IEEEproof}

Figure~\ref{ParametricPb} seems to illustrate the presence of a minimum for $P_b$ as a function of $s$. But this result is unreachable and conjectured here. Nevertheless, we propose the following weak result.
\begin{Proposition}\label{extremum=1}
For any $\lambda/\mu>0$, $P_b$ has an extremum $P_b\text{*}$  for $\rho_V=1$. 
\end{Proposition}

\begin{IEEEproof}
Existence is clear since $P_b$ is continuous and, by definition, $P_b(0)=1$ and  $P_b$ tends to $1$ as $\rho_V$ tends to $ +\infty$ using equation~\eqref{symP_b}.

The value $\rho_V=1$ comes also from the symmetry proved in Property~\ref{Prop:symmetry}. Indeed, using equation ~\eqref{symP_b},
\begin{align}
P_b'(\rho_V) &= -\frac{1}{\rho_V^2} P_b'(1/\rho_V) \text{.} \nonumber
\end{align}
Therefore, $P_b'(1)=-P_b'(1)=0$ and the  sign of the derivative changes at $\rho_V=1$. Hence, for $\rho_V=1$ a local extremum is reached for $P_b$, and the same holds for $P_b^{+}$.
\end{IEEEproof}
Uniqueness of the extremum should come from the  convexity of the parametric curve $\rho_V \mapsto (s(\rho_V), P_b(\rho_V))$. This yields that it is a minimum. Nevertheless this convexity result is still to be proved because the implicit relation between $\rho_V$ and $\rho_R$ (see Section~\ref{sec:SteadyStateBehaviour}) makes calculations tedious.

\begin{Conj}
\label{Th:minimum=1}
For any $\lambda/\mu>0$, $P_b$ has a unique minimum $P_b\text{*}$ reached for $\rho_V=1$.
\end{Conj}

Conjecture~\ref{Th:minimum=1} is similar to the result for  bike-sharing systems (see ~\cite{Fricker} for details), that the minimum is reached in $\rho_V=1$. In Figure~\ref{ParametricPb}, we observe the existence and uniqueness of the minimum conjectured in~\ref{Th:minimum=1}. The curves plotted for both systems are convex with a unique minimum.

Moreover, for the influence of the fleet size on  the system behavior, the fact that $P_b$ tends to $1$ as $\rho_V$ tends to $ +\infty$ makes an important difference with the homogeneous bike-sharing system. Without reservation, $P_b$ goes to $1$ when $s$ tends to infinity whereas reservation implies full saturation when $s=K$. This fact will be developed in the following with the study of the heavy-traffic case.

\begin{Prop}
\label{Prop:PbRhoVfixed}
For all $\rho_V > 0$, the approximated proportion of problematic stations $P_b^{+}(\rho_V)$ is increasing with traffic load $\lambda/\mu$.
\end{Prop}
\begin{IEEEproof}
Define  $a=\lambda/\mu > 0$  and $f_a=f/a$  where $f$ is given by equation~\eqref{f}. It holds that
\begin{align}\label{f_a}
f_a(\rho_V, \rho_R) = (1-\rho_R/a) Z(\rho_V, \rho_R) - \sum\limits_{l=0}^K \frac{\rho_R^l}{l!}.
\end{align}
Denote by $\psi_a$ the diffeomorphism $\psi$ defined in Theorem~\ref{Theorem:diffeomorphism}. 

Let $\alpha > 1$ be fixed. First let us  prove that 
\begin{align}\label{ineq}
\alpha \psi_a > \psi_{\alpha a}.
\end{align} 
It is equivalent to prove that $f_{\alpha a}(\rho_V, \alpha \rho_R) < 0$ where $(\rho_V, \rho_R)$ different to $(0,0)$ satisfies $f_a(\rho_V, \rho_R) =0$, i.e. $\psi_a(\rho_V) = \rho_R$. For that, using equation~\eqref{f_a}, it holds that
\begin{multline*}Z(\rho_V, \rho_R) f_{\alpha a}(\rho_V, \alpha \rho_R) = \\
Z(\rho_V, \rho_R) (1 - \rho_R/a) \sum\limits_{(k,l) \in \chi} \rho_V^k \frac{\rho_R^l}{l!}\alpha^l - Z(\rho_V, \rho_R) \sum\limits_{l=0}^K \frac{\rho_R^l}{l!} \alpha^l.
\end{multline*}

Since $f_a(\rho_V, \rho_R)=0$, $Z(\rho_V, \rho_R) (1 - \rho_R/a) = \sum\limits_{l=0}^K \frac{\rho_R^l}{l!}$. Hence,
\begin{align}
&  Z(\rho_V, \rho_R) f_{\alpha a}(\rho_V, \alpha \rho_R) \nonumber\\
				& =  \sum\limits_{l=0}^K \frac{\rho_R^l}{l!} \sum\limits_{(k,l) \in \chi} \rho_V^k \frac{\rho_R^l}{l!}\alpha^l  -  \sum\limits_{(k,l) \in \chi} \rho_V^k \frac{\rho_R^l}{l!} \sum\limits_{l=0}^K \frac{\rho_R^l}{l!} \alpha^l \nonumber \\
				& =  \sum\limits_{k=0}^K \sum\limits_{l=0}^{K-k} \sum\limits_{j=0}^{K-k} \rho_V^k \frac{\rho_R^{j+l}}{j!l!} ( \alpha^l - \alpha^j ) \nonumber\\
				& + \sum\limits_{k=0}^K \sum\limits_{l=0}^{K-k} \sum\limits_{j=K-k+1}^K \rho_V^k \frac{\rho_R^{j+l}}{j!l!} ( \alpha^l - \alpha^j ). \label{NegSum}
\end{align}
The first term of the right-hand side of the equation~\eqref{NegSum} equals zero as it is symmetrical in $l$ and $j$. And because $j > l$ and $\alpha > 1$ the second term of the right-hand side of equation~\eqref{NegSum} is negative, which concludes the first part of this proof.

Then, by equation~\eqref{eq:Pbexpression}, 
\begin{align*}
P_b^{+}(\alpha a) = 2 - \frac{\psi_{\alpha a}(\rho_V)}{\alpha a} - \frac{\psi_{\alpha a}(\rho_V)}{\alpha a \rho_V}.
\end{align*}
 Plugging equation~\eqref{ineq} in it, we can conclude that
$P_b^{+}(\alpha a) > P_b^{+}(a)$. It ends the proof.
\end{IEEEproof}

Property~\ref{Prop:PbRhoVfixed} is not so interesting in practice since the operator does not have directly access to the parameter $\rho_V$. However, it could be a first step for proving the following Conjecture~\ref{Conj:Pbsfixed}.

\begin{Conj}
\label{Conj:Pbsfixed}
For any  $s \in \mathopen{[} 0, K \mathclose{[}$, the proportion of problematic stations $P_b^{+}(s)$ is increasing with traffic $\lambda/\mu$.
\end{Conj}

We conjecture both Property~\ref{Prop:PbRhoVfixed} and Conjecture~\ref{Conj:Pbsfixed} are valid for $P_b$ as well. In fact, assertions given in Conjectures~\ref{Conj:Pbsfixed} and~\ref{Th:minimum=1} for $P_b$ can be observed on Figure~\ref{ParametricPb} where $P_b$ is plotted as a function of $s$.

Conjecture~\ref{Conj:Pbsfixed} would imply the behavior of the system with reservation is very different from systems without reservation. Indeed, for bike-sharing systems, increasing traffic does not systematically worsens the situation (\textit{e.g.} if $s$ is close to $K$ as we can see on  Figure~\ref{ParametricPb} (see also~\cite{Fricker}). But in both systems, increasing the capacity improves the system performance.

In \cite{Fricker}, for the corresponding bike-sharing model, the minimum $2/(K+1)$ for $P_b$ is reached for $\rho_V=K/2+\lambda/\mu$. We are also interested   analytical expressions for this minimum for the car-sharing system  but it turns out to be tedious again because of implicit equations. Therefore, even if we would have proven that this minimum is reached for $\rho_V=1$, an asymptotic expression could only be obtained in two cases: light traffic ($\lambda/\mu \rightarrow 0$) or heavy traffic ($\lambda/\mu \rightarrow +\infty$).

\begin{Prop}[Light traffic]
\label{Prop:DASaPetitToutRhoV}
For any $\rho_V \geq 0$, when $\lambda/\mu$ tends to  $0$,
\begin{align}
 s &= \sum\limits_{k=0}^K k\rho_V^k/\sum\limits_{k=0}^K \rho_V^k + o(1) \nonumber\\
P_b &=  1 - \frac{\rho_V^K-\rho_V}{\rho_V^{K+1}-1} + o(1) \nonumber
\end{align}
\end{Prop}
\begin{IEEEproof}
The method  is standard. It uses that $\rho_R = \frac{1-\rho_V^K}{1-\rho_V^{K+1}}\rho_V \frac{\lambda}{\mu} + o(\lambda/\mu)$. The proof is omitted.
\end{IEEEproof}

Property~\ref{Prop:DASaPetitToutRhoV} means that, in case of light traffic ($\lambda/\mu \rightarrow 0$), the system has the same performance at first order in $\lambda/\mu$ as the homogeneous bike-sharing model 
. This result coincides with the intuition that when the traffic is low, reservation does not have time to induce congestion because cars are picked-up faster than new users arrive in the system. Particularly, the optimal limiting proportion of problematic stations is obtained when $\rho_V = 1$ which gives that $P_b^{*} = 2/(K+1)$ and $s^{*}=K/2$ at first order in $\lambda/\mu$. This result can be observed on Figure~\ref{ParametricPb} where curves for the bike-sharing system get quite close to curves for the car-sharing system when the traffic rate  decreases.

\begin{Prop}[Optimality for light traffic]
\label{Prop:DASaPetitRhov=1}
If $\rho_V=1$, as $\lambda/\mu \rightarrow 0$,
\begin{align*}
s &= \frac{K}{2}+ \frac{K^2}{2(K+1)^2}  \frac{\lambda}{\mu}  + O((\lambda/\mu)^2 ) , \\
P_b &= \frac{2}{K+1} + \frac{2K}{(K+1)^3} \frac{\lambda}{\mu}  + O((\lambda/\mu)^2 ).
\end{align*}
Assuming Conjecture~\ref{Th:minimum=1}, the previous quantities are respectively  $s^*$ and $P_b^*$.
\end{Prop}
\begin{IEEEproof}
These formulas are straightforwardly obtained, via $\rho_R = K/(K+1)\; \lambda/\mu+ O((\lambda/\mu)^2 )$.
\end{IEEEproof}

Assuming Conjecture~\ref{Th:minimum=1}, Property~\ref{Prop:DASaPetitRhov=1} highlights the differences of behavior between the bike-sharing system and the system with reservation at the optimum. Indeed, $s^{*}$ increases less rapidly in the present case and, most of all, $P_b^{*}$ is no more constant but has a positive term in $\lambda/\mu$ when traffic is light. However, increasing the capacity $K$ dramatically reduces the effects of this term.

\begin{Prop}[Heavy traffic]
\label{Prop:DASaGrandToutRhoV}
For all $\epsilon >0$ and $M>\epsilon$, for any $\rho_V \in [ \epsilon, M]$, as $\lambda/\mu \rightarrow +\infty$,
\begin{align*}
 s &= K - \sqrt{\frac{K}{\rho_V}}\sqrt{\mu/\lambda}  + \frac{3+\rho_V-2K}{2\rho_V}\mu/\lambda \\
		&+ O (\mu/\lambda)^{3/2}) \\
P_b &=  1 - (K-1)\left(\frac{\lambda}{\mu} \right)^{-1} + \frac{(\rho_V+1)(K-1)}{\sqrt{\rho_V K}}\left(\frac{\lambda}{\mu} \right)^{-3/2} \\
		&+ O \left(\left(\frac{\lambda}{\mu} \right)^{-2} \right) 
\end{align*}
\end{Prop}
\begin{IEEEproof}
We proof fisrt that $\rho_R = \sqrt{K\rho_V} \sqrt{\lambda/\mu}  -  (1+\rho_V)/2+ \frac{(5/8-K/2)(1+\rho_V^2)+\rho_V/4}{\sqrt{\rho_V K}}\sqrt{\mu/\lambda}  + O (\mu/\lambda)$. The method  is standard. 
\end{IEEEproof}

The behaviour of $P_b$ is interesting in Property~\ref{Prop:DASaGrandToutRhoV} as it happens to tend to 1 rather quickly compared to the convergence of $s$ and independently of $\rho_V$ (the first two terms do not depend on $\rho_V$). This is also the reason why we chose to dig further in calculations than in Property~\ref{Prop:DASaPetitToutRhoV}. We shall also mention that this speed of convergence is higher because of the term $\rho_R^K/(K!Z)$ which does not exist in $P_b^{+}$. Assuming Conjecture~\ref{Th:minimum=1} to be true, it is interesting to consider the case $\rho_V=1$. 

\begin{Prop}[Optimality for heavy traffic]
\label{Prop:DASaGrandRhov=1}
If $\rho_V=1$, as $\lambda/\mu \rightarrow +\infty$,
\begin{align*}
s &=  K - \sqrt{\frac{K}{\lambda/\mu}}+ \frac{2-K}{\lambda/\mu} +  O \left(\left(\frac{\lambda}{\mu} \right)^{-3/2} \right) \\
P_b &= 1 - \frac{K-1}{\lambda/\mu}  + \frac{2(K-1)}{\sqrt{K}}\left(\frac{\lambda}{\mu} \right)^{-3/2} + O \left(\left(\frac{\lambda}{\mu} \right)^{-2} \right).
\end{align*}
Assuming Conjecture~\ref{Th:minimum=1}, the previous quantities are respectively  $s^*$ and $P_b^*$.
\end{Prop}
\begin{IEEEproof}
We obtain these formulas straightforwardly from Property~\ref{Prop:DASaGrandToutRhoV}.
\end{IEEEproof}

\begin{figure}[h]
\centering
\begin{tikzpicture}[scale=0.8]
\begin{axis}[
title={Problematic stations at a given $s$, with $K=4$.},
area legend,
legend style={at={(0.97,0.03)},
anchor=south east},
axis x line=bottom,
axis y line=left,
xlabel={Traffic $\lambda/\mu$},
ylabel={Proportion of problematic stations $P_b$},
xmin=0,
ymin=0,
xmax=10,
ymax=1
]

\legend{$s\text{*}(\lambda/\mu)$,$s=1.5$, $s=2.1$, $s=2.6$, $s=3.25$}
\addplot[black, dashed, line legend] table {./ProblematicStations=fct_a_rhoV=1_K=4.txt};
\addplot[black, line legend] table {./ProblematicStations=fct_a_s0=1.5_K=4.txt};
\addplot[blue, line legend] table {./ProblematicStations=fct_a_s0=2.1_K=4.txt};
\addplot[red, line legend] table {./ProblematicStations=fct_a_s0=2.6_K=4.txt};
\addplot[green, line legend] table {./ProblematicStations=fct_a_s0=3.25_K=4.txt};
\end{axis}
\end{tikzpicture}
\caption{Evolution of the proportion of problematic stations with traffic ($s$ fixed). Solid lines represent this evolution for different values of $s$. The dashed line represents the evolution if for each $\lambda/\mu$ the optimal parameter $s^{*}$ is chosen. Consistently with Figure~\ref{ParametricPb}, the lowest curve depends on the traffic.}
\label{fig:ProblematicStationsSfixed}
\end{figure}
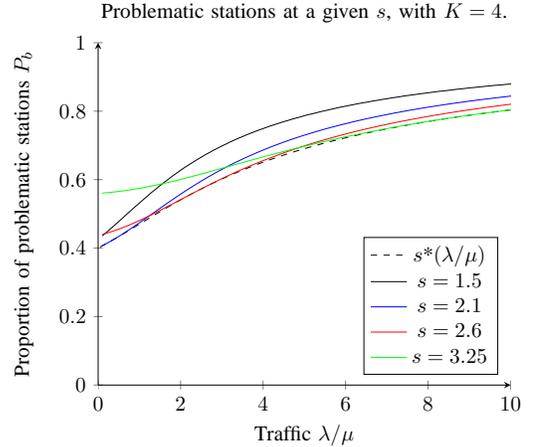

\subsection{Maximizing the number of trips}

The operator of the system, either a private company or the city council, may choose parameters differently since it has to take into account the cost of running the system. In this case, the number of successful trips  becomes more relevant as a metric because the revenue of the operator is more or less directly proportional.  While it is just an approximated value, we assume in the following that $T = \lambda (1-P_b)$.

 Let us assume $\mu$ is fixed. We will study $T$ as a function of $\lambda$ which is common in transportation.
Let us first study the maximum number of trips $T^{*}=\lambda(1-P_b^*)$ that  can be obtained as a function of $\lambda$. For that, the previous expansions will be very useful.
We use for $\lambda$ small the results for light traffic of Property~\ref{Prop:DASaPetitRhov=1}, which gives
\begin{align*}
T^*=\lambda \left(1-\frac{2}{K+1}\right)+O(\lambda^2)
\end{align*}
which means a linear function of the demand for small $\lambda$.
For heavy traffic case, from  Property~\ref{Prop:DASaGrandRhov=1}, as $\lambda$ tends to $\infty$,
\begin{align*}
T^*=(K-1)\mu+\frac{2(K-1)\mu^{3/2}}{\sqrt{K\lambda}}+O\left(\frac{1}{\sqrt{\lambda}}\right)
\end{align*}

\begin{figure}[h]
\centering
\begin{tikzpicture}[scale=0.8]
\begin{axis}[
title={Successful trips at a given $s$, with $K=4$.},
area legend,
legend style={at={(0.97,0.03)},
anchor=south east},
axis x line=bottom,
axis y line=left,
xlabel={Traffic $\lambda/\mu$},
ylabel={Successful Trips $T$},
xmin=0,
ymin=0,
xmax=10,
ymax=2.4
]

\legend{$s\text{*}(\lambda/\mu)$,$s=1.5$, $s=2.1$, $s=2.6$, $s=3.25$}
\addplot[black, dashed, line legend] table {./SuccessfulTrips=fct_a_rhoV=1_K=4.txt};
\addplot[black, line legend] table {./SuccessfulTrips=fct_a_s0=1.5_K=4.txt};
\addplot[blue, line legend] table {./SuccessfulTrips=fct_a_s0=2.1_K=4.txt};
\addplot[red, line legend] table {./SuccessfulTrips=fct_a_s0=2.6_K=4.txt};
\addplot[green, line legend] table {./SuccessfulTrips=fct_a_s0=3.25_K=4.txt};
\end{axis}
\end{tikzpicture}
\caption{Number of successful trips $T$ as a function of traffic $\lambda/\mu$ for different values of fixed $s$ (solid lines) and for the optimal value of $s$ (dashed line).}
\label{fig:SuccessfulTrips}
\end{figure}
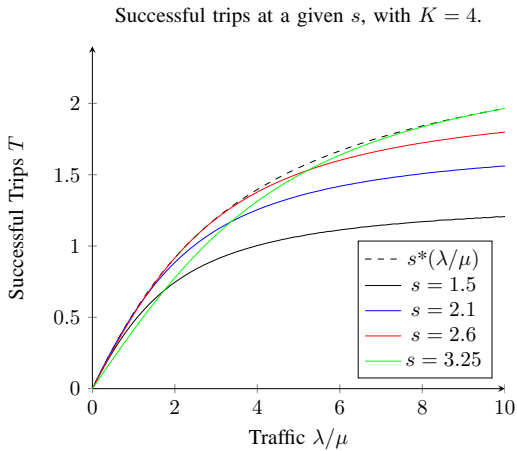

Figure~\ref{fig:SuccessfulTrips} shows the evolution of the number of successful trips as a function of $\lambda/\mu$. In practice, it is obtained by solving numerically equations~\eqref{eq:RhoR} and~\eqref{eq:sMoyenne} for each value of $\lambda/\mu$. We observe the linearity when $\lambda/\mu \rightarrow 0$. When $\lambda/\mu \rightarrow \infty$, if $s$ is fixed, $\rho_V \rightarrow 0$ which prevents from using Property~\ref{Prop:DASaGrandToutRhoV} ($\rho_V$ has to be bounded). This is consistent with the fact the solid curves do not seem to tend to $K-1$ whereas the dashed curve does (since it corresponds to $\rho_V=1$ fixed). Another interesting observation is that there is no local extremum ($T$ is strictly increasing) which is an argument in favour of Conjecture~\ref{Conj:Pbsfixed}.

Demand $\lambda$ is a function of the price $p$ of a trip. The profit is

$$\Pi(s,p) = \lambda(p) (1-P_b(s,p)) - C(s) \text{,}$$

where $C$ is the cost function. Maximizing profit leads to choose $s$ such that 

$$\frac{\partial P_b}{\partial s}= - \frac{1}{\lambda(p)} C'(s) \text{.}$$

Since $C$ is supposed to be strictly increasing in $s$ and Conjecture~\ref{Th:minimum=1} to hold, the optimal choice in terms of profit for the operator is necessarily reached for some $s_0 < s^{*}$. Thus, even under competition pressure, the trade-off between price and number of bikes certainly results in choosing $s < s^{*}$.

\section{Extension to a Model with Double Reservation}
\label{sec:DoubleReservation}

Car-sharing systems often offer the opportunity to reserve online a car and a parking space in the desired station, a while before actually picking-up the car. In this section we present a model that adresses this particular demand of travel.

\subsection{Model Description}


At a given station $i$, new reservations of car are made at rate $\lambda$. At the same time, the user making this reservation also wants to reserve a parking space in some station $j$. If there is no car available in $i$ \textit{or} no available parking space in $j$, the user is rejected. Otherwise, she waits a time exponentially distributed with mean $1/\nu$ before coming to pick-up her car. Then, the journey towards station $j$ takes an exponentially distributed time with mean $1/\mu$. The user returns her car at station $j$ and leaves the system.

A three-dimensional state space (reserved cars, reserved spaces and free cars) is unsuited to fully describe this process as markovian because the duration of reservation of spaces is the sum of two exponentially distributed variables (the time to pick-up the car and the time to travel). This leads to introducing a fourth variable to distinguish between spaces reserved by users \textit{not yet travelling} and users \textit{travelling}. As we will see in the description of transitions, the Markov process associated to this model is quite complicated. Let us make an approximation, introduce the following process and discuss later the relevance of the approximation.

Let us denote now $\chi= \{ (j,k,l,m) \in \mathbb{N}^4, j+k+l+m \leq K \}$ and $Y^N_{j,k,l,m}(t)$ the proportion of the $N$ stations where, at time $t$: 
\begin{center}
\begin{tabular}{rclp{4cm}}
$V_i^{r,N}(t)$ & $=$ & $j$ & Reserved vehicles \\
$V_i^N(t)    $ & $=$ & $k$ & Free vehicles \\
$R^N(t)      $ & $=$ & $l$ & Reserved parking spaces by users travelling\\
$R_i^{r,N}(t)$ & $=$ & $m$ & Reserved parking spaces by users not yet travelling
\end{tabular}
\end{center}

which can be written
$$Y^N_{j,k,l,m}(t) = \frac{1}{N} \sum\limits_{i=1}^N \textbf{1}_{\{V_i^{r,N}(t)=j;V_i^N(t)=k;R_i^N(t)=l;R_i^{r,N}(t)=m \}} \text{.}$$

As the system is homogeneous, the process $(Y^N(t)) = (Y^N_{j,k,l,m}(t))_{(j,k,l,m) \in \chi}$ is a Markov process. Suppose the process is at state $(y_{j,k,l,m})_{(j,k,l,m) \in \chi}$. There are three different types of transitions:

\begin{itemize}
\item \textbf{Reservation}. At rate $\lambda$, a user wants to reserve a car in a station $i$. The station is in state $(j,k,l,m)$ with probability $N y_{j,k,l,m}$. Simultaneously, she chooses an arrival station $i'$ in state $(j',k',l',m')$ with probability $y_{j',k',l',m'}$. She succeeds if there is at least one free car in $i$ ($k>0$) and at least one free space in $i'$ ($j'+k'+l'+m'<K$). Therefore, the transition rate is $\lambda N y_{j,k,l,m} y_{j',k',l',m'}$ if $k>0$ and $j'+k'+l'+m'<K$. The car reservation causes $y_{j,k,l,m}$ to decrease and $y_{j+1,k-1,l,m}$ to increase, both by $1/N$. The parking reservation causes $y_{j',k',l',m'}$ to decrease and $y_{j',k',l',m'+1}$ to increase, both by $1/N$.

\item \textbf{Cars picked up}. Time between reservation and picking up is exponentially distributed with parameter $\nu$. For stations in state $(j,k,l,m)$, there are $j N y_{j,k,l,m}$ reserved cars. Hence, cars are picked up at rate $j \nu y_{j,k,l,m}N$. But, at the same time, a reserved parking space changes of status (from user \textit{not yet travelling} to \textit{travelling}) in another \textit{random} station somewhere else, sampled with probability proportional to the number of reserved parking spaces in this station. Thus, it occurs in a station in state $(j',k',l',m')$ with probability

$$ \frac{m' N y_{j',k',l',m'}}{\sum\limits m' N y_{j',k',l',m'}} = y_{j',k',l',m'} \frac{m'}{\mathbb{E}_y (R^{r,N})}$$

where  $\mathbb{E}_y [R^{r,N}]$ is the average number of reserved parking space for users not yet travelling. Finally, the transition rate is $\nu j \frac{m'}{\mathbb{E}_y [R^{r,N}]} N y_{j,k,l,m}y_{j',k',l',m'}$. In the original model, this change of status occurs in the very same station as the one chosen during reservation: the approximation allows to keep the memoryless property. The transition causes $y_{j,k,l,m}$ to decrease and $y_{j-1,k,l,m}$ to increase. The change of status in the arrival station causes $y_{j',k',l',m'}$ to decrease and $y_{j',k',l'+1,m'-1}$ to increase. All by $1/N$.

\item \textbf{Cars returned}. Considering stations in state $(j,k,l,m)$, the number of reserved parking spaces is $l N y_{j,k,l}$.  As trip are exponentially distributed with mean $1/\mu$, cars are returned at rate $\mu l N y_{j,k,l}$. When a car arrives, it causes $y_{j,k,l,m}$ to decrease and $y_{j,k+1,l-1,m}$ to increase, both by $1/N$.
\end{itemize}


\subsection{Probabilistic Interpretation}
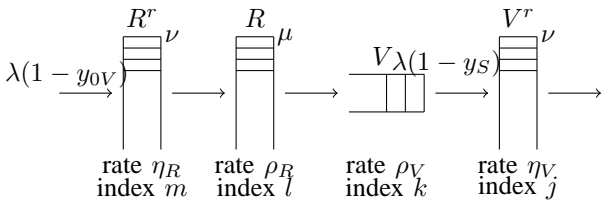
\begin{figure}[h]
\centering
\begin{tikzpicture}[scale=0.5]

\node (arrival) at (0,0.5) {$\lambda (1-y_{0V})$};
\draw (arrival);
\draw[->] (-0.2,0) -- (1.2,0);

\draw[-] (1.5,-1.5) -- (1.5,1.5);
\draw[-] (2.5,-1.5) -- (2.5,1.5);
\draw[-] (1.5,1.5) -- (2.5, 1.5);
\draw[-] (1.5,1.2) -- (2.5, 1.2);
\draw[-] (1.5,0.9) -- (2.5, 0.9);
\draw[-] (1.5,0.6) -- (2.5, 0.6);

\node (Rr) at (2,2) {$R^r$};
\draw (Rr);
\node (etaR) at (2,-2) {rate $\eta_R$};
\draw (etaR);
\node (IndexetaR) at (2,-2.5) {index $m$};
\draw (IndexetaR);


\node (middle1) at (2.8,1.5) {$\nu$};
\draw (middle1);
\draw[->] (2.8,0) -- (4.2,0);

\draw[-] (4.5,-1.5) -- (4.5,1.5);
\draw[-] (5.5,-1.5) -- (5.5,1.5);
\draw[-] (4.5,1.5) -- (5.5, 1.5);
\draw[-] (4.5,1.2) -- (5.5, 1.2);
\draw[-] (4.5,0.9) -- (5.5, 0.9);
\draw[-] (4.5,0.6) -- (5.5, 0.6);

\node (R) at (5,2) {$R$};
\draw (R);
\node (roR) at (5,-2) {rate $\rho_R$};
\draw (roR);
\node (IndexroR) at (5,-2.5) {index $l$};
\draw (IndexroR);


\node (middle2) at (5.8,1.5) {$\mu$};
\draw (middle2);
\draw[->] (5.8,0) -- (7.2,0);

\draw[-] (7.5,.5) -- (9.5,.5);
\draw[-] (7.5,-.5)-- (9.5,-.5);
\draw[-] (9.5,.5) -- (9.5, -.5);
\draw[-] (9,.5) -- (9, -.5);
\draw[-] (8.5,.5) -- (8.5, -.5);

\node (V) at (8.4,1) {$V$};
\draw (V);
\node (roV) at (8.5,-2) {rate $\rho_V$};
\draw (roV);
\node (IndexroV) at (8.5,-2.5) {index $k$};
\draw (IndexroV);


\node (middle3) at (10.1,0.8) {$\lambda (1-y_S)$};
\draw (middle3);
\draw[->] (9.8,0) -- (11.2,0);

\draw[-] (11.5,-1.5) -- (11.5,1.5);
\draw[-] (12.5,-1.5) -- (12.5,1.5);
\draw[-] (11.5,1.5) -- (12.5, 1.5);
\draw[-] (11.5,1.2) -- (12.5, 1.2);
\draw[-] (11.5,0.9) -- (12.5, 0.9);
\draw[-] (11.5,0.6) -- (12.5, 0.6);

\node (Vr) at (12,2) {$V^r$};
\draw (Vr);
\node (etaR) at (12,-2) {rate $\eta_V$};
\draw (etaR);
\node (IndexetaV) at (12,-2.5) {index $j$};
\draw (IndexetaV);


\node (out) at (12.8,1.5) {$\nu$};
\draw (out);
\draw[->] (12.8,0) -- (14.2,0);

\end{tikzpicture}
\caption{A typical station as a tandem of 4 queues with overall capacity $K$.}
\label{fig:MMqueuesModel3}
\end{figure}

The simplification introduced in the second transition allows to write down a differential system of equations similar to what is done for the simple-reservation model: the limiting empirical distribution $y(t)$ of the stations evolves in time as the distribution of some non-homogeneous Markov process on $\chi$, whose jumps are given by $L_{y(t)}$, updated by the current distribution $y(t)$. This process can be seen as the queue length vector of four coupled queues with an overall capacity $K$ (see Figure~\ref{fig:MMqueuesModel3}).

As it is done for the simple-reservation model, this interpretation leads to 4 equations on the parameters of the queues ($\rho_R$,$\rho_V$,$\eta_R$,$\eta_V$), and a fifth equation is obtained from the expression of the average number of vehicles in the system $s(\rho_R,\rho_V,\eta_R,\eta_V) = s$. Simple algebra can reduce most of them to the equations studied in the former model via a change of variables. Indeed, first, $\eta_R$ and $\eta_V$ can be omitted as they appear to be directly proportional to $\rho_R$. Remain 3 equations on $\rho_R$ and $\rho_V$. Second, if $a=\lambda/\mu$, 

\begin{align*}
\tilde{\rho}_R & = \left(1+2\frac{\mu}{\nu}\right)\rho_R \\
\text{and} \qquad \tilde{a} & = \left(1+2\frac{\mu}{\nu}\right) a
\end{align*}

leads to $\tilde{\rho}_R$ and $\rho_V$ satisfying the exact same equations \ref{eq:RhoR} and \ref{eq:RhoV}, with parameter $\tilde{a}$. However, regarding the relation with the average number of vehicles, the result is not as simple since after the change of variable we have

\begin{equation}
\label{eq:sDoubleReservation}
s(\rho_V, \rho_R) = \frac{1}{Z} \sum\limits_{k+l \leq K} \left(k+l\frac{1+\mu/\nu}{1+2\mu/\nu}\right) \rho_V^k \frac{\tilde{\rho}_R^l}{l!}.
\end{equation}

The increase of $s$ with respect to its parameters still has to be proven because of this new coefficient $\left(k+l\frac{1+\mu/\nu}{1+2\mu/\nu}\right)$. But if we assume $1/\nu$ to be sufficiently small, an argument of continuity allows to use Theorem~\ref{Theorem:sCroissante}. Therefore, existence and mostly uniqueness of the equilibrium point are straightforward with the different results we have for the simple-reservation model. 

\subsection{Consequences on Performance}
Thanks to the change of variable, we also have that $P_b$ writes the same way as a function of $\rho_V$ and $\tilde{\rho}_R$, with parameter $\tilde{a}$. As a consequence, we still have that a local minimum of $P_b$ is reached for $\rho_V = 1$ as mentioned in Section~\ref{sec:Performance} and the value of this minimum is the same as for the first model with a traffic of $\tilde{a}= \left(1+2\frac{\mu}{\nu}\right)$. Similarly, $P_b(0)=1$ and $P_b \rightarrow 1$ when $\rho_V \rightarrow +\infty$.

However, Equation~\ref{eq:sDoubleReservation} shows that for some number of cars $s$, the associated solutions $\tilde{\rho}_R$ and $\rho_V$ have to be higher than the corresponding $\rho_R$ and $\rho_V$ in the first model. Therefore, if we compare $s \mapsto P_b(s)$ to its analogue in the first model with traffic $\tilde{a}$, it is smaller for low values of $s$, reaches a minimum before and then is higher for $s$ closer to $K$.

\subsection{Relevancy of the model}
This model introduced above is simplified in the description of the transition of users picking up their cars. However, given the probabilistic interpretation in terms of four coupled queues, we can intuitively expect it to have the same stationary behavior as the non-simplified process. To confirm this, we performed simulations of the original model where people actually go to their reserved parking space instead of re-sampling one, and derived the steady-state distribution. 
The empirical performance function resulting from this latter distribution $s \mapsto P_b(s)$ seems identical to the one given by the model.

\section{Conclusion and Future Work}
\label{sec:Conclusion}
In this paper we have investigated the influence of reservation on the performance of vehicle-sharing systems. We used a stochastic model and then a mean-field limit to study the steady-state behaviour of the system  and to compare it with former homogeneous models without reservation. Though reservation makes calculations more tedious, assuming convergence, we still have theoretical results on the existence and uniqueness  of the limiting stationary state of one station when the system gets large. Techniques are standard with a crucial use of probabilistic tools in each step. Because results on limiting performance, in particular the optimal performance, are difficult to obtain, we gave results in the two cases of light and heavy traffic. It appears to be the same as for bike-sharing systems when traffic is light. But, unlike in bike-sharing systems, saturation occurs when traffic is high or when the average number of bikes per station $s$ is close to the capacity. We conjectured and partially proved that the optimal choice of $s$ is reached for $\rho_V=1$, like for bikes, which is close  to $K/2+\lambda/(2\mu)$ in the light traffic case, when $\lambda/\mu$ is very small, but no more in general. In heavy traffic, when $\lambda/\mu$ is large, it is of the order of $K-\sqrt{\frac{K}{\lambda/\mu}}$. 
\\\indent
To take into account real possibilities in such systems, we also studied a network where users can reserve both a car and a space before beginning their journey. We studied an approximated model (validated by simulations) and we proved this new model is mathematically similar to the previous model: all is as if the traffic was increased.

Some conjectures presented in this paper still have to be proved. And similarly to what has been done for bike-sharing systems in \cite{Fricker}, further work might focus on different algorithms to improve performance or to generalize the model. For instance, in terms of the so called passive regulation, we investigated the power of two choice in bike-sharing systems and it could also have an important impact on car-sharing systems. Or, concerning the regulation by the operator, since trucks cannot be used with cars we could look at different strategies for \textit{ambassadors}: they move cars \textit{one by one} from full stations to empty ones which would make the model easier to define than for trucks that move several bikes together. In terms of generalization, for bike-sharing systems  (see~\cite{Heterogeneous}), results extend to an inhomogeneous model  with clusters of stations with the same parameters. The same can be done for our model.

\bigskip


\begin{thebibliography}{10}

\bibitem{Carlier-1}
A.~Carlier, A.~Munier Kordon, and Witold Klaudel.
\newblock Optimization of a one-way carsharing system with relocation
  operations.
\newblock In {\em Proceedings of 10th International Conference on Modeling,
  Optimization and Simulation}, Nancy, France, November 5-7 2014.

\bibitem{Chemla}
Daniel Chemla, Frédéric Meunier, and Roberto~Wolfler Calvo.
\newblock Bike sharing systems: Solving the static rebalancing problem.
\newblock {\em Discrete Optimization}, 10(2):120 -- 146, 2013.

\bibitem{Fayolle-7}
G.~Fayolle and J.-M. Lasgouttes.
\newblock Asymptotics and scalings for large product-form networks via the
  central limit theorem.
\newblock {\em Markov Process. Related Fields}, 2(2):317--348, 1996.

\bibitem{Fricker}
Christine Fricker and Nicolas Gast.
\newblock Incentives and redistribution in homogeneous bike-sharing systems
  with stations of finite capacity.
\newblock {\em European Journal on Transportation and Logistics (accepted)},
  2014.

\bibitem{Heterogeneous}
Christine Fricker, Nicolas Gast, and Hanene Mohamed.
\newblock Mean field analysis for inhomogeneous bike sharing systems.
\newblock {\em DMTCS Proceedings, 23rd Intern. Meeting on Probabilistic,
  Combinatorial, and Asymptotic Methods for the Analysis of Algorithms
  (AofA'12)}, 2012.

\bibitem{Tibi}
Christine Fricker and Danielle Tibi.
\newblock Equivalence of ensembles for large vehicle-sharing models.
\newblock Preprint, august 2014.

\bibitem{Simatos}
Ayalvadi Ganesh, Sarah Lilienthal, D.~Manjunath, Alexandre Proutiere, and
  Florian Simatos.
\newblock Load balancing via random local search in closed and open systems.
\newblock In {\em Proceeding SIGMETRICS '10 Proceedings of the ACM SIGMETRICS
  international conference on Measurement and modeling of computer systems},
  pages 287--298, 2010.

\bibitem{gast2010}
N.~Gast and B.~Gaujal.
\newblock {A mean field model of work stealing in large-scale systems}.
\newblock {\em ACM SIGMETRICS Performance Evaluation Review}, 38(1):13--24,
  2010.

\bibitem{George-1}
David~K. George and Cathy~H. Xia.
\newblock Asymptotic analysis of closed queueing networks and its implications
  to achievable service levels.
\newblock {\em SIGMETRICS Performance Evaluation Review}, 38(2):3--5, 2010.

\bibitem{Kaspi-1}
Mor Kaspi, Tal Raviv, Michal Tzur, and Hila Galili.
\newblock Regulating vehicle sharing systems through parking reservation
  policies: Analysis and performance bounds.
\newblock working paper, august 2014.

\bibitem{Robert}
Philippe Robert.
\newblock {\em Stochastic Networks and Queues}.
\newblock Stochastic Modelling and Applied Probability Series. Springer-Verlag,
  New York, 2003.

\bibitem{Jost}
Ariel Waserhole, Vincent Jost, and Nadia Brauner.
\newblock Pricing techniques for self regulation in vehicle sharing systems.
\newblock {\em Electronic Notes in Discrete Mathematics}, 41(0):149 -- 156,
  2013.

\bibitem{Bogenberger-1}
S.~Weikl and K.~Bogenberger.
\newblock Relocation strategies and algorithms for free-floating car sharing
  systems.
\newblock In {\em 15th International IEEE Conference on Intelligent
  Transportation Systems}, Anchorage, Alaska, USA, September 16-19 2012.

\end{thebibliography}

\bibliographystyle{IEEEtran}

\section{Appendices}
\label{sec:Appendices}

\subsection{Proof of Property~\ref{sCroissante}}
\begin{IEEEproof}
We will prove that, for all $(\rho_V, \rho_R) \in \Gamma$, 
$$\frac{\partial \tilde{s}}{\partial \rho_R}(\rho_V, \rho_R) >0 \text{ and }  \frac{\partial \tilde{s}}{\partial \rho_V}(\rho_V, \rho_R) >0$$
 by induction on $K$. For that, let $\tilde{s}$ be denoted by $\tilde{s}_K$. 

By a change of indexes, $\tilde{s}_K$ can be rewritten
$$\tilde{s}_K  = \frac{\sum\limits_{k+l \leq K} (k+l)\rho_V^k \frac{\rho_R^l}{l!}}{\sum\limits_{k+l \leq K} \rho_V^k \frac{\rho_R^l}{l!}} =  \frac{\sum\limits_{k=0}^K k p_k}{\sum\limits_{k=0}^K p_k} $$
where, by definition,
$$p_k = \sum\limits_{i=0}^k \rho_V^i \frac{\rho_R^{k-i}}{(k-i)!} \text{.}$$
Define also for $(k,l)\in \chi$
$$r_{l,k} = \frac{p_l}{p_k}.$$

Let $k>0$ be fixed. We first show that $r_{k,k-1}$ is an increasing function of both $\rho_R$ and $\rho_V$. Indeed

$$r_{k,k-1} = \frac{p_k}{p_{k-1}} = \frac{ \rho_V p_{k-1} + \rho_R^{k}/k! }{p_{k-1}} = \rho_V + \frac{\rho_R^{k}}{k!}\frac{1}{p_{k-1}} \text{,}$$

thus
\begin{align}\label{dRhoV}
\frac{\partial r_{k,k-1}}{\partial \rho_V} = 1 - \frac{\rho_R^{k}}{k!} \frac{\partial p_{k-1}}{\partial \rho_V} \frac{1}{p_{k-1}^2}.
\end{align}

But, by definition of $p_k$,
\begin{align}\label{arche}
\rho_R^{k} \frac{\partial p_{k-1}}{\partial \rho_V} = \sum\limits_{i=1}^{k-1} i \rho_V^{i-1}\frac{\rho_R^{2k-i-1}}{(k-i-1)!} = \sum\limits_{i=0}^{k-2} (i+1) \rho_V^{i}\frac{\rho_R^{2k-i-2}}{(k-i-2)!}.
\end{align}

And, using that all terms of the sum in the following equation are positive,

\begin{align}\label{biche}
k!p_{k-1}^2 &= k! \sum\limits_{1 \leq u,v \leq k} \rho_V^{u+v-2} \frac{\rho_R^{2k-u-v}}{(k-u)!(k-v)!} \nonumber\\
	&> \sum\limits_{i=0}^{k-2} \rho_V^{i}\rho_R^{2k-i-2} \sum\limits_{j=1}^{i+1} \frac{k!}{(k-2-i+j)!(k-j)!} \text{.} 
\end{align}

For all $j,\;0\leq j  \leq i+1$, while
 $k-2-i+j < k$, it holds that

$$\frac{k!}{(k-2-i+j)!(k-j)!} > \frac{1}{(k-i-2)!}$$
and then 
$$\sum\limits_{j=1}^{i+1} \frac{k!}{(k-2-i+j)!(k-j)!} > \frac{i+1}{(k-i-2)!} \text{.}$$
Plugging in equation~\eqref{biche} and comparing with  equation~\eqref{arche}, it gives
$$k!p_{k-1}^2 > \rho_R^{k} \frac{\partial p_{k-1}}{\partial \rho_V}.$$ 
Therefore, using ~\eqref{dRhoV}, it allows to conclude that $$\frac{\partial r_{k,k-1}}{\partial \rho_V} >0 \text{.}$$

Moreover,
\begin{align}\label{chat}
\frac{\partial r_{k,k-1}}{\partial \rho_R} = \frac{\rho_R^{k-1}}{(k-1)!p_{k-1}} \left( 1 - \frac{\rho_R}{kp_{k-1}}\frac{\partial p_{k-1}}{\partial \rho_R} \right) \text{.}
\end{align}
Using that
$$\frac{\partial p_{k-1}}{\partial \rho_R}=\sum\limits_{i=0}^{k-1}\rho_V^{i}\frac{\rho_R^{k-i-2}}{(k-1-i)!}(k-1-i)= p_{k-2}$$ 
and $kp_{k-1} > \rho_R p_{k-2}$, because
\begin{align*}
\frac{\rho_R}{k} p_{k-2}&=\sum\limits_{i=0}^{k-2}\rho_V^{i}\frac{\rho_R^{k-i-1}}{(k-2-i)!k}\\
&<\sum\limits_{i=0}^{k-2}\rho_V^{i}\frac{\rho_R^{k-i-1}}{(k-1-i)!}\\ 
&<p_{k-1},
\end{align*}
we can conclude that $\frac{\partial r_{k,k-1}}{\partial \rho_R} >0$.

Consequently, if $l>k$, $r_{l,k} = \prod\limits_{i=k+1}^l r_{i,i-1}$ is an increasing function of $\rho_V$ and $\rho_R$. This gives that $u_K$ defined by 

$$u_K=\frac{p_K}{\sum\limits_{k=0}^K p_k}=\frac{1}{\sum\limits_{k=0} r_{k,K}}$$
is strictly increasing in $x \in \{ \rho_V, \rho_R \}$, because $r_{k,K}=1/r_{K,k}$ is strictly decreasing in $x$.

Let us mention that $\tilde{s}_0$ is constant with $x$ and that
$$\tilde{s}_K = (1-u_K) \tilde{s}_{K-1} + K u_K,$$
which yields that
$$\frac{\partial \tilde{s}_K}{\partial x} = (K-\tilde{s}_{K-1}) \frac{\partial u_K}{\partial x} + (1-u_K) \frac{\partial \tilde{s}_{K-1}}{\partial x} \text{.}$$

Since $K-\tilde{s}_{K-1} > 0$, $u_K <1$ and $\frac{\partial u_K}{\partial x} >0$, by induction we can conclude that $\frac{\partial \tilde{s}_K}{\partial x}>0$ for all $K \geq 1$. It ends the proof.
\end{IEEEproof}

\end{document}